\title{On integral and differential formulations in nonlocal elasticity}

\date{}

\documentclass[12pt]{article}

\usepackage{amsmath}
\usepackage{amssymb}
\usepackage{xcolor}

\usepackage{graphicx}
\usepackage{float}
\usepackage{overpic}

\author{J.Kaplunov, D.A. Prikazchikov and  L. Prikazchikova$^*$
\\{\small \it School of Computing and Mathematics,}\\
{\it \small Keel University, Keele, ST5 5BG, UK}\\
{\small $^*$Corresponding author: l.prikazchikova@keele.ac.uk}}

\begin{document}

\maketitle

\begin{abstract}

The paper is concerned with comparative analysis of differential and integral formulations for boundary value problems in nonlocal elasticity. For the sake of simplicity, the focus is on an antiplane problem for a half-space with prescribed shear stress along the surface. In addition,  1D exponential kernel depending on the vertical coordinate is considered.

First, a surface loading in the form of a travelling harmonic wave is studied. This provides a counter-example, revealing that within the framework of Eringen's theory the solution to the differential model does not satisfy the equation of motion in nonlocal stresses underlying the related integral formulation.

A more general differential setup, starting from singularly perturbed equations expressing the local stresses through the nonlocal ones, is also investigated. It is emphasized that the transformation of the original integral formulation to the differential one in question is only possible provided that two additional conditions on nonlocal stresses hold on the surface. As a result, the formulated problem subject to three boundary conditions appears to be ill-posed, in line with earlier observations for equilibrium of a nonlocal cantilever beam.

Next, the asymptotic solution of the  singularly perturbed problem, subject to a prescribed stress on the boundary, together with only one of the aforementioned extra conditions, is obtained at a small internal size. Such simplification may be justified when only one of the stresses demonstrates nonlocal behaviour; a similar assumption has been recently made within the so-called dilatational gradient elasticity. Three-term expansion is obtained, leading to a boundary value problem in local stresses over interior domain. The associated differential equations are identical to those proposed by Eringen, however, the derived effective boundary condition incorporates the effect of a nonlocal boundary layer which has previously been  ignored. 
Moreover, the calculated nonlocal correction to the classical antiplane problem for an elastic half-space, coming from the boundary conditions is by order of magnitude greater than that appearing in the equations of motion. Finally, it is shown that the proposed effective condition supports an antiplane surface wave.

\end{abstract}

\section{Introduction}
Nonlocal elasticity is an important theory finding various modern applications including nano-technology, e.g. see  \cite{elishakoff2013applications, karlicic2015nonlocal, peddieson2003application} and references therein. The concept of nonlocal stresses and integral relations expressing nonlocal stresses through their local counterparts, have been introduced by A.C.Eringen, see the most influential contributions \cite{eringen1972nonlocal, eringen1983differential}. In the paper \cite{eringen1983differential} the original integral nonlocal theory has been  reduced to a differential form, which is much easier for implementation. The latter has been adapted for numerous scenarios, in particular for thin nano-structures, see \cite{reddy2007nonlocalplates, hache2019asymptotic, aksencer2011levy} to name a few. In addition, the differential formulation in nonlocal elasticity formally has a lot of in common with a popular model of gradient elasticity, e.g. see \cite{askes2011gradient} and references therein.

Nowadays, there is a strong belief in the equivalence of integral and differential nonlocal setups among a broad academic community. At the same time, the highly cited paper \cite{eringen1983differential} as well as later publications on the subject do not verify whether the solutions of nonlocal differential equations satisfy the original integral relations. 


This paper is aimed at bridging  the gap in understanding of the relationship between integral and differential approaches in nonlocal elasticity. For the sake of simplicity, antiplane problem for a nonlocally elastic half-space is considered for a particular one-dimensional exponential kernel, see \cite{eringen1983differential}. We begin with the example of travelling harmonic waves, induced by a prescribed surface stress. Elementary calculations readily indicate that the solution of the conventional differential model, introduced in  \cite{eringen1983differential}, does not satisfy the equation of motion for nonlocal stresses underlying the integral theory. 

Next, we consider a more general differential model, starting from singularly perturbed differential equations, expressing classical (local) stresses through nonlocal ones. In this case, the coefficients at senior derivatives corresponding to an internal size are assumed to be small in comparison with a typical wavelength. The considered boundary value problem involves two boundary conditions along the surface of the half-space, in contrast to the aforementioned differential nonlocal theory \cite{eringen1983differential}, as well as the classical local antiplane elasticity, both operating with a single boundary condition. At the same time, it may be easily shown that the associated differential equations are equivalent to the initial integral relations only provided that each of the nonlocal stresses satisfies an extra boundary condition at the surface. Thus, in total we obtain three boundary conditions. This observation is in line with the conclusion of \cite{romano2017constitutivebc}, demonstrating that an elastostatic problem in nonlocal integral elasticity for a cantilever beam is ill-posed, see also more recent publications \cite{mikhasev2020solution, vaccaro2021limit}.
Below, we proceed taking into account only one of two extra boundary conditions, evaluating a posteriori the discrepancy arising from the violation of the remaining boundary condition.

A three-term asymptotic solution is derived. Along with a slowly varying component, it includes a boundary layer localised near the surface. 
Similar boundary layers arise in integral nonlocal elasticity, e.g. see 
\cite{chebakov2016nonlocalhalfspace, chebakov2017nonlocalplate, abdollahi2014nonlocal}, but are usually ignored within the differential formulation, e.g. see \cite{karlicic2015nonlocal}, dealing with nonlocally elastic engineering structures. 
At the same time, boundary layers have been taken into consideration in gradient elasticity, see \cite{lurie2021dilatation}.

Asymptotic analysis results in effective boundary conditions to the equations in antiplane elasticity, expressed in terms of local stresses. As in previous asymptotic developments within integral nonlocal elasticity for a Gaussian kernel \cite{chebakov2016nonlocalhalfspace, chebakov2017nonlocalplate}, a leading order correction to the conventional (local) boundary conditions is due to the effect of boundary layers. This correction is of order of magnitude greater than that coming from the nonlocal equations of motions.

As an example, the derived effective boundary conditions are implemented for analysis of antiplane shear surface waves. It is worth noting that the differential nonlocal model \cite{eringen1983differential} does not support such waves, which have been also observed in recent publications using different methodologies, e.g. see \cite{eremeyev2020strongly}.

\section{Antiplane problem in nonlocal elasticity for a half space}

\begin{figure}[H]
\centering
\begin{overpic}[scale=0.9]{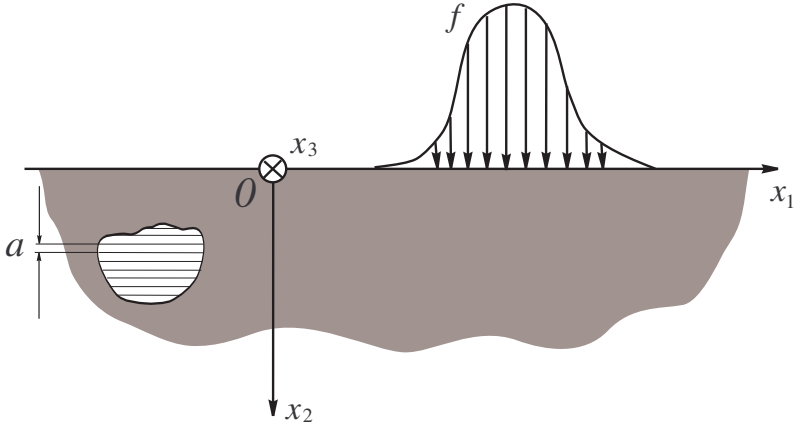}
\end{overpic}
\caption{Antiplane shear of a nonlocally elastic half-space. }
\label{fig_nonlocal}
\end{figure}

Consider antiplane shear for an elastic half-space ($-\infty < x_1, x_3 < +\infty$, $0 \leq x_2 < +\infty$), see Figure \ref{fig_nonlocal}. 

The governing equation of motion is written in terms of nonlocal stresses as
\begin{equation}\label{eq_mot}
\dfrac{\partial s_{13}}{\partial x_1}+\dfrac{\partial s_{23}}{\partial x_2}=\rho \dfrac{\partial^2 u}{\partial t^2},
\end{equation}
subject to the boundary condition
\begin{equation}\label{eq_bc}
s_{23}=f(x_1,t) \quad \mbox{at} \quad x_2=0.
\end{equation}
Here $s_{j3}=s_{j3}(x_1,x_2,t)$, $j=1,2$ are nonlocal stresses, $u=u(x_1,x_2,t)$ is out of plane displacement, $\rho$ is mass density and $t$ is time.
For the sake of simplicity, the nonlocal stresses are taken in the form
%
%
%
\begin{equation}\label{eq_s1D}
s_{j3}=\dfrac{1}{2a}\int_{0}^{\infty}e^{-\frac{|x_2-x_2'|}{a}} \sigma_{j3}(x_1, x_2',t) \operatorname{d}x_2', \quad j=1,2,
\end{equation}
adapting 1D kernel from \cite{eringen1983differential} (cf. formula [3.4]). In the above, $a$ is a small parameter (e.g. lattice parameter) and local stresses $\sigma_{j3}$ are expressed through displacement as
\begin{equation}\label{eq_sigma}
\sigma_{j3}=\mu \dfrac{\partial u}{\partial x_j},
\end{equation}
with $\mu$ being the Lam\'e shear modulus.

As shown in \cite{eringen1983differential} for an infinite media, equation of motion  \eqref{eq_mot} can be transformed into the following differential form
\begin{equation}\label{eq_motd}
\dfrac{\partial \sigma_{13}}{\partial x_1}
+\dfrac{\partial \sigma_{23}}{\partial x_2}
=\rho\left(\dfrac{\partial^2 u}{\partial t^2}-a^2 
\dfrac{\partial^4 u}{\partial t^2 \partial x_2^2}\right).
\end{equation}
In the current paper we aim to compare the integral and differential formulations for the bounded domain. We begin with a counter-example, demonstrating that for a half-space the analysed differential and integral formulations are generally not equivalent.
%
%
\section{Counter-example}
Let us consider a simple example of a forced time-harmonic problem for antiplane motion of a half-space, motivated by analogous treatment of surface waves in \cite{eringen1983differential}.
The load along $x_2=0$ is taken as
\begin{equation}
f(x_1,t)=F e^{ik_0(x_1-c_0 t)},
\end{equation}
where amplitude $F$, wavenumber $k_0$ and speed $c_0$ are prescribed quantities, such that $c_0<c_2$ and $k_0 \ll 1/a$, with $c_2 = \sqrt{\mu / \rho}$ denoting the shear wave speed.

Below we operate with the displacement $u$ as a function of transverse variable only $u=u(x_2)$, assuming that the exponential factor $e^{ik_0(x_1-c_0 t)}$ has already been separated. Hence, equation \eqref{eq_motd} becomes
\begin{equation}\label{eq_motu}
\left(c_2^2-a^2 k_0^2 c_0^2\right)\dfrac{\operatorname{d}^2 u}{\operatorname{d} x_2^2}-k_0^2 \left(c_2^2 -c_0^2\right) u=0,
\end{equation}
Boundary condition \eqref{eq_bc} may be re-written as
\begin{equation}\label{eq_bcu}
\int_{0}^{\infty}e^{-\frac{|x_2-x_2'|}{a}} \dfrac{\operatorname{d}u}{\operatorname{d} x_2'} \operatorname{d}x_2' = \dfrac{2aF}{\mu} \quad \mbox{at} \quad x_2=0.
\end{equation}
The decaying solution of \eqref{eq_motu} is 
\begin{equation}\label{eq_exp}
u(x_2) = U e^{-k_0\gamma x_2},
\end{equation}
where $U$ is an arbitrary constant and
\begin{equation}\label{eq_gamma}
\gamma = \sqrt{\dfrac{c_2^2-c_0^2}{c_2^2-a^2k_0^2c_0^2}}.
\end{equation}
Using boundary condition \eqref{eq_bcu}, the sought for solution for the displacement is found as 
\begin{equation}\label{eq_u2}
u(x_2) = -\dfrac{2F(1+ak_0\gamma)}{k_0\mu\gamma}e^{-k_0\gamma x_2}.
\end{equation}
The limit $\gamma \rightarrow 0$ when displacement $u \rightarrow \infty$ does not correspond to a decaying solution of interest, similarly to the classical antiplane problem. Using \eqref{eq_s1D}, the nonlocal stresses may be expressed as
\begin{align}
&s_{13} = \dfrac{iF}{\gamma(1-ak_0\gamma)}\left\{(1+ak_0 \gamma)e^{-\frac{x_2}{a}}-2e^{-k_0\gamma x_2} \right\}, \label{eq_s12}\\
&s_{23} = -\dfrac{F}{(1-ak_0\gamma)}\left\{(1+ak_0\gamma)e^{-\frac{x_2}{a}}-2e^{-k_0\gamma x_2} \right\}\label{eq_s32}.
\end{align}
Now we need to verify whether solutions \eqref{eq_u2}, \eqref{eq_s12} and \eqref{eq_s32} satisfy the original equation of motion in nonlocal stresses. Multiplying them  by $e^{ik_0(x_1-c_0 t)}$,   we deduce from \eqref{eq_mot}
\begin{equation}
e^{-\frac{x_2}{a}}\dfrac{F(1+ak\gamma)(\gamma-ak)}{a\gamma(1-ak\gamma)}  \not\equiv 0.
\end{equation}
Hence, this equation does not hold, up to a boundary layer of width $\sim a$, except for a rather special loading for which $\gamma = ak_0$. Thus, it has been demonstrated that the solutions \eqref{eq_u2}-\eqref{eq_s32} of the differential model generally do not satisfy  the equation of motion \eqref{eq_mot} within the original framework of integral nonlocal elasticity.

\section{Singularly perturbed differential model}
Let us now approach the issue of equivalence of differential and integral formulations in nonlocal elasticity from a more general perspective. Following the methodology introduced in  \cite{eringen1983differential}, we differentiate relations \eqref{eq_s1D} twice with respect to $x_2$, obtaining  singularly perturbed equations
\begin{equation}\label{eq_sigmas_diff}
a^2 \dfrac{\partial^2 s_{j3}}{\partial x_2^2}-s_{j3} = -\sigma_{j3}, \quad j=1,2.
\end{equation}
Expressing the local stresses $\sigma_{j3}$ through the nonlocal ones and substituting the result into \eqref{eq_s1D}, we derive additional conditions on the boundary
\begin{equation}\label{eq_addcond}
\left.\left(s_{j3}-a \dfrac{\partial s_{j3}}{\partial x_2}\right)\right|_{x_2=0}=0, \quad j=1,2.
\end{equation}
In view of these additional conditions  equivalence of two approaches is hardly possible. Indeed, there is already a condition \eqref{eq_bc} on the boundary, in addition to two extra conditions \eqref{eq_addcond}. For the moment let us satisfy one of the conditions \eqref{eq_addcond}, say the one for ($j=1$), and evaluate at later stage the discrepancy caused by the violation of the remaining condition ($j=2$).

In contrast to a traditional differential model in Section 3, we consider three differential equations \eqref{eq_mot} and\eqref{eq_sigmas_diff}, and two boundary conditions \eqref{eq_bc} and \eqref{eq_addcond} for $j=1$, instead of the equation of motion \eqref{eq_motd} in local stresses subject to boundary condition \eqref{eq_bc} only. For both cases local stresses are expressed through displacements by \eqref{eq_sigma}.

In fact, on differentiating \eqref{eq_sigmas_diff} twice with respect to time and using equation of motion \eqref{eq_mot} with relation \eqref{eq_sigma} it is possible to obtain two coupled equations in nonlocal stresses, namely
\begin{align}
\begin{split}
& \dfrac{\partial^2 s_{13}}{\partial x_1^2}+\dfrac{\partial^2 s_{23}}{\partial x_1 \partial x_2}-\dfrac{1}{c_2^2}  \dfrac{\partial^2}{\partial t^2} \left(  s_{13}-a^2 \dfrac{\partial^2 s_{13}}{\partial x_2^2}\right) =0,\\
&\dfrac{\partial^2 s_{23}}{\partial x_2^2}+\dfrac{\partial^2 s_{13}}{\partial x_1 \partial x_2}-\dfrac{1}{c_2^2}  \dfrac{\partial^2}{\partial t^2} \left(  s_{23}-a^2 \dfrac{\partial^2 s_{23}}{\partial x_2^2}\right) =0.
\end{split}
\end{align}
These second order equations (in spacial variables) obviously require two boundary conditions along the surface.
However, for the sake of clarity, in what follows we operate with singularly perturbed equations \eqref{eq_sigmas_diff} and equation of motion \eqref{eq_mot}.


%
%
\section{Asymptotic analysis}

The formulated problem can be tackled  asymptotically due to the presence of a natural small parameter
\begin{equation}
\eta = \dfrac{a}{l} \ll 1,
\end{equation}
where $l$ is a typical wavelength. 
Below, dimensionless variables
\begin{equation}
\zeta_p=\frac{x_2}{l}, \quad \zeta_q=\frac{x_2}{a}, \quad \xi_1=\frac{x_1}{l}, \quad \tau=t\frac{c_2}{l}, 
\end{equation}
are introduced, with slow and fast variables $\zeta_p$ and $\zeta_q$, associated with two types of variation of the unknown quantities over the transverse coordinate. 
Also, the dimensionless functions
\begin{equation}
u^{*}=\frac{u}{l}, \quad \sigma_{j3}^{*}=\frac{\sigma_{j3}}{\mu}, \quad s_{j3}^{*}=\frac{s_{j3}}{\mu}, \quad f^*=\dfrac{f}{\mu}, \quad j=1,2,
\end{equation}
are chosen, assuming all the starred quantities to be of order unity. 

Now, the governing equations may be rewritten as
\begin{align}
&\dfrac{\partial s_{13}^*}{\partial \xi_1}  +\dfrac{\partial s_{23}^*}{\partial \zeta_p}= \dfrac{\partial^2 u^*}{\partial \tau^2},\\
&\eta^2 \dfrac{\partial^2 s_{j3}^*}{\partial \zeta_p^2} -s_{j3}^* = -\sigma_{j3}^*, \quad j=1,2,
\end{align}
subject to the boundary conditions
\begin{align}
&s_{23}^*\big|_{\zeta_p=0}=f^*(\xi_1,\tau),\\
&s_{13}^*\big|_{\zeta_p=0}=\eta\left. \dfrac{\partial s_{13}^*}{\partial \zeta_p} \right|_{\zeta_p=0}.
\end{align}
Similarly to \cite{chebakov2017nonlocalplate}, the nonlocal stresses are split into fast and slow components $p_{j3}^*$ and $q_{j3}^*$, respectively, i.e.  
\begin{equation}
s_{13}^* = p_{13}^* + q_{13}^*
\end{equation}
and
\begin{equation}\label{eq_s23_pq}
s_{23}^* = p_{23}^* + \eta q_{23}^*,
\end{equation}
where $p_{ij}^*=p_{ij}^*(\xi_1,\zeta_p, \tau)$ and $q_{ij}^*=q_{ij}^*(\xi_1,\zeta_q, \tau)$. 
Hence, the equations above are re-cast in terms of fast and slow quantities as
\begin{align}
&\dfrac{\partial p_{13}^*}{\partial \xi_1}+\dfrac{\partial p_{23}^*}{\partial \zeta_p}  =   \dfrac{\partial^2 u^*}{\partial \tau^2},\\
&\dfrac{\partial q_{13}^*}{\partial \xi_1}+\dfrac{\partial q_{23}^*}{\partial \zeta_q}=0, \\
&p_{j3}^*-\eta^2 \dfrac{\partial^2 p_{j3}^*}{\partial \zeta_p^2} = \sigma_{j3}^*,\\
&q_{j3}^*- \dfrac{\partial^2 q_{j3}^*}{\partial \zeta_q^2} =0, \quad j=1,2,
\end{align}
subject to
\begin{align}
&p_{13}^*\big|_{\zeta_p=0} + q_{13}^*\big|_{\zeta_q=0} = \left. \dfrac{\partial q_{13}^*}{\partial \zeta_q}\right|_{\zeta_q=0}
+\eta \left.\dfrac{\partial p_{13}^*}{\partial \zeta_p} \right|_{\zeta_p=0}\\
&p_{23}^*\big|_{\zeta_p=0} + \eta q_{23}^*\big|_{\zeta_q=0} = f^*.
\end{align}

At leading (zero) order we have
\begin{align}
&\dfrac{\partial p_{13}^{(0)}}{\partial \xi_1}+\dfrac{\partial p_{23}^{(0)}}{\partial \zeta_p}  =   \dfrac{\partial^2 u^{(0)}}{\partial \tau^2},\label{eq_lo1}\\
&\dfrac{\partial q_{13}^{(0)}}{\partial \xi_1}+\dfrac{\partial q_{23}^{(0)}}{\partial \zeta_q}=0, \label{eq_lo2}\\
&p_{j3}^{(0)} = \sigma_{j3}^{(0)}\label{eq_lo3},\\
&q_{j3}^{(0)}- \dfrac{\partial^2 q_{j3}^{(0)}}{\partial \zeta_q^2} =0, \quad j=1,2,\label{eq_lo4}
\end{align}
together with
\begin{align}
&p_{13}^{(0)}\big|_{\zeta_p=0} + q_{13}^{(0)}\big|_{\zeta_q=0} = \left.\dfrac{\partial q_{13}^{(0)}}{\partial \zeta_q}\right|_{\zeta_q=0},\label{eq_lo5}\\
&p_{23}^{(0)}\big|_{\zeta_p=0} = f^*. \label{eq_lo6}
\end{align}
Combining \eqref{eq_lo1} and \eqref{eq_lo3}, we deduce
\begin{equation}\label{eq_losigma}
\dfrac{\partial \sigma_{13}^{(0)}}{\partial \xi_1}+\dfrac{\partial \sigma_{23}^{(0)}}{\partial \zeta_p}=\dfrac{\partial^2 u^{(0)}}{\partial \tau^2},
\end{equation}
subject to 
\begin{equation}
\sigma_{23}^{(0)} \big|_{\zeta_p=0}=f^*.
\end{equation}
In additon, stress-displacement relations at each order follow from \eqref{eq_sigma} as
\begin{equation}\label{eq_sigmau_k}
\sigma_{13}^{(k)}= \dfrac{\partial u^{(k)}}{\partial \xi_1}  \quad \mbox{and} \quad \sigma_{23}^{(k)}=\dfrac{\partial u^{(k)}}{\partial \zeta_p}, \quad k=0,1,2,\dots
\end{equation}
Now, from \eqref{eq_lo4} it follows that
\begin{equation}
q_{j3}^{(0)}=Q_j^{(0)}(\xi_1,\tau)e^{-\zeta_q}.
\end{equation}
Substituting the latter into \eqref{eq_lo2}, we get
\begin{equation}
Q_2^{(0)} = \dfrac{\partial Q_1^{(0)}}{\partial \xi_1}.
\end{equation}
Next, using condition \eqref{eq_lo5},
\begin{equation}
Q_{1}^{(0)} = -\dfrac{1}{2}\sigma_{13}^{(0)} \big|_{\zeta_p=0}.
\end{equation}
Therefore, in terms of local stresses
\begin{align}
&q_{13}^{(0)} = -\dfrac{1}{2}\sigma_{13}^{(0)} \big|_{\zeta_p=0} \, e^{-\zeta_q},\\
&q_{23}^{(0)} = -\dfrac{1}{2} \left. \dfrac{\partial\sigma_{13}^{(0)}}{\partial \xi_1}\right|_{\zeta_p=0} \, e^{-\zeta_q}.
\end{align}

At first order, we have
\begin{align}
&\dfrac{\partial p_{13}^{(1)}}{\partial \xi_1}+\dfrac{\partial p_{23}^{(1)}}{\partial \zeta_p}  =   \dfrac{\partial^2 u^{(1)}}{\partial \tau^2},\label{eq_no1}\\
&\dfrac{\partial q_{13}^{(1)}}{\partial \xi_1}+\dfrac{\partial q_{23}^{(1)}}{\partial \zeta_q}=0, \label{eq_no2}\\
&p_{j3}^{(1)} = \sigma_{j3}^{(1)}\label{eq_no3},\\
&q_{j3}^{(1)}- \dfrac{\partial^2 q_{j3}^{(1)}}{\partial \zeta_q^2} =0,  \quad j=1,2, \label{eq_no4}
\end{align}
subject to the following boundary conditions
\begin{align}
&p_{13}^{(1)}\big|_{\zeta_p=0} + q_{13}^{(1)}\big|_{\zeta_q=0} = \dfrac{\partial p_{13}^{(0)}}{\partial \zeta_p}\bigg|_{\zeta_p=0} +
\dfrac{\partial q_{13}^{(1)}}{\partial \zeta_q}\bigg|_{\zeta_q=0}, \label{eq_no5}\\
&p_{23}^{(1)}\big|_{\zeta_p=0} + q_{23}^{(0)}\big|_{\zeta_q=0} = 0. \label{eq_no6}
\end{align}
Using \eqref{eq_no3} together with \eqref{eq_no1} we obtain
\begin{equation}\label{eq_no_u}
\dfrac{\partial \sigma_{13}^{(1)}}{\partial \xi_1}+\dfrac{\partial \sigma_{23}^{(1)}}{\partial \zeta_p}=\dfrac{\partial^2 u^{(1)}}{\partial \tau^2},
\end{equation}
and
\begin{equation}
\sigma_{23}^{(1)} \big|_{\zeta_p=0}=\dfrac{1}{2}\dfrac{\partial \sigma_{13}^{(0)}}{\partial \xi_1}\bigg|_{\zeta_p=0}.
\end{equation}
Also, we have as before from \eqref{eq_no4} and \eqref{eq_no2},
\begin{equation}
q_{j3}^{(1)}=Q_j^{(1)}(\xi_1,\tau)e^{-\zeta_q} \quad \mbox{and} \quad
Q_2^{(1)} = \dfrac{\partial Q_1^{(1)}}{\partial \xi_1}.
\end{equation}
Next, using boundary condition \eqref{eq_no5}, we deduce
\begin{equation}
Q_{1}^{(1)} = \dfrac{1}{2}\left(
\dfrac{\partial \sigma_{13}^{(0)}}{\partial \zeta_p} \bigg|_{\zeta_p=0}-\sigma_{13}^{(1)} \big|_{\zeta_p=0}\right).
\end{equation}
Thus, all the first order terms are determined. 

At next order we have
\begin{align}
&\dfrac{\partial p_{13}^{(2)}}{\partial \xi_1}+\dfrac{\partial p_{23}^{(2)}}{\partial \zeta_p}  =   \dfrac{\partial^2 u^{(2)}}{\partial \tau^2},\label{eq_nno1}\\
&\dfrac{\partial q_{13}^{(2)}}{\partial \xi_1}+\dfrac{\partial q_{23}^{(2)}}{\partial \zeta_q}=0, \label{eq_nno2}\\
&p_{j3}^{(2)} - \dfrac{\partial^2 p_{j3}^{(0)}}{\partial \zeta_p^2} = \sigma_{j3}^{(2)}\label{eq_nno3},\\
&q_{j3}^{(2)}- \dfrac{\partial^2 q_{j3}^{(2)}}{\partial \zeta_q^2} =0, \quad j=1,2,\label{eq_nno4}
\end{align}
together with
\begin{align}
&p_{13}^{(2)}\big|_{\zeta_p=0} + q_{13}^{(2)}\big|_{\zeta_q=0} = \dfrac{\partial q_{13}^{(2)}}{\partial \zeta_q} \bigg|_{\zeta_q=0}+\dfrac{\partial p_{13}^{(1)}}{\partial \zeta_p}\bigg|_{\zeta_p=0},\label{eq_nno5}\\
&p_{23}^{(2)}\big|_{\zeta_p=0} + q_{23}^{(1)}\big|_{\zeta_q=0} = 0. \label{eq_nno6}
\end{align}
Substituting \eqref{eq_nno3} into \eqref{eq_nno1} and taking into account \eqref{eq_losigma}, we obtain
\begin{equation}\label{eq_nno_u}
\dfrac{\partial \sigma_{13}^{(2)}}{\partial \xi_1}+\dfrac{\partial \sigma_{23}^{(2)}}{\partial \zeta_p}=\dfrac{\partial^2 u^{(2)}}{\partial \tau^2}-\dfrac{\partial^4 u^{(0)}}{\partial \zeta_p^2 \partial \tau^2}.
\end{equation}
This equation is solved subject to
\begin{equation}
\sigma_{23}^{(2)} \big|_{\zeta_p=0}=-
\dfrac{\partial^2 \sigma_{23}^{(0)}}{\partial \zeta_p^2}\bigg|_{\zeta_p=0} -\dfrac{1}{2}
\dfrac{\partial^2 \sigma_{13}^{(0)}}{\partial \xi_1 \partial \zeta_p}\bigg|_{\zeta_p=0}+\dfrac{1}{2}\dfrac{\partial \sigma_{13}^{(1)}}{\partial \xi_1} \bigg|_{\zeta_p=0}.
\end{equation}
As before, it follows from   \eqref{eq_nno4}, \eqref{eq_nno2} and \eqref{eq_no5} that
\begin{equation}
q_{j3}^{(2)}=Q_j^{(2)}(\xi_1,\tau)e^{-\zeta_q},
\end{equation}
where
\begin{equation}
Q_{1}^{(2)} = \dfrac{1}{2}\left(-\dfrac{ \partial^2 \sigma_{13}^{(0)}}{\partial \zeta_p^2}\bigg|_{\zeta_p=0}+
\dfrac{\partial \sigma_{13}^{(1)}}{\partial \zeta_p} \bigg|_{\zeta_p=0} - \sigma_{13}^{(2)} \big|_{\zeta_p=0}\right)
\end{equation}
and
\begin{equation}
Q_2^{(2)} = \dfrac{\partial Q_1^{(2)}}{\partial\xi_1}.
\end{equation}
%
%
%

A combination of three equations, corresponding to leading, first and second orders, see \eqref{eq_losigma}, \eqref{eq_no_u} and \eqref{eq_nno_u}, implies
\begin{equation}
\dfrac{\partial \sigma_{13}^*}{\partial \xi_1}+\dfrac{\partial \sigma_{23}^*}{\partial \zeta_p}=\dfrac{\partial^2 u^*}{\partial \tau^2}-\eta^2 \dfrac{\partial^4 u^*}{\partial \zeta_p^2 \partial \tau^2},
\end{equation}
subject to boundary condition at $\zeta_p=0$
\begin{equation}
\sigma_{23}^* - \dfrac{\eta}{2} \dfrac{ \partial \sigma_{13}^*}{\partial \xi_1} + \eta^2 \left(\dfrac{1}{2}\dfrac{\partial^2\sigma_{13}^*}{\partial \xi_1 \partial \zeta_p} + \dfrac{\partial^2 \sigma_{23}^*}{\partial \zeta_p^2}  \right)= f^*. 
\end{equation}
In the above
\begin{align}
&\sigma_{j3}^*=\sigma_{j3}^{(0)}+\eta\sigma_{j3}^{(1)}+\eta^2 \sigma_{j3}^{(2)},\\
&u^*=u^{(0)}+\eta u^{(1)}+\eta^2 u^{(2)}.
\end{align}

\section{Discussion}
In terms of the original variables these equations can be re-cast as
\begin{equation}\label{eq_motfin}
\dfrac{\partial \sigma_{13}}{\partial x_1}+\dfrac{\partial \sigma_{23}}{\partial x_2}=\rho \left(\dfrac{\partial^2 u}{\partial t^2}-a^2 \dfrac{\partial^4 u}{\partial t^2 \partial x_2^2 }\right),
\end{equation}
with
\begin{equation}\label{eq_bcfin}
\sigma_{23} - \dfrac{a}{2} \dfrac{ \partial \sigma_{13}}{\partial x_1} + a^2 \left(\dfrac{1}{2}\dfrac{\partial^2\sigma_{13}}{\partial x_1 \partial x_2} + \dfrac{\partial^2 \sigma_{23}}{\partial x_2^2}  \right)= f, \end{equation}
at $x_2=0$, where stresses are expressed through displacements via \eqref{eq_sigma}. 

Thus, an original problem of nonlocal elasticity has been transformed to that expressed in terms of local stresses. Over a narrow near-surface zone of width of order $\eta=a/l$, the solution to the latter problem should be complemented by the boundary layer component, calculated in Section 5. The equation of motion  \eqref{eq_motfin} was initially derived in \cite{eringen1983differential}. However,  the $O(\eta)$ correction to the boundary conditions is of order of magnitude higher than $O(\eta^2)$ correction to the equations of motion.
%
%
%

Let us now return back to condition  \eqref{eq_addcond} for $j=2$, which has not been taken into account throughout the asymptotic derivation in Section 5. 
Retaining the terms up to $O(\eta^2)$, we have from \eqref{eq_s23_pq} and subsequent formulae in Section 5
\begin{align}
\begin{split}
s_{23}^*&=p_{23}^{(0)}+\eta\left(p_{23}^{(1)}+q_{23}^{(0)} \right)+\eta^2\left(p_{23}^{(2)}+q_{23}^{(1)}  \right)\\
&=\sigma_{23}^{(0)}+\eta\left(\sigma_{23}^{(1)}-\dfrac{1}{2} \dfrac{\partial \sigma_{13}^{(0)}}{\partial \xi_1}\bigg|_{\zeta_p=0}e^{-\zeta_q} \right) \\ 
&\qquad +\eta^2\left(\sigma_{23}^{(2)}+\dfrac{\partial^2\sigma_{23}^{(0)}}{\partial \zeta_p^2}+\dfrac{1}{2}\left(\dfrac{\partial^2 \sigma_{13}^{(0)}}{\partial \xi_1 \partial \zeta_p}\bigg|_{\zeta_p=0}-\dfrac{\partial\sigma_{13}^{(1)}}{\partial \xi_1}\bigg|_{\zeta_p=0} \right)e^{-\zeta_q} \right),
\end{split}
\end{align}
or, equivalently
\begin{equation}\label{eq_s32diff}
s_{23}^*=\sigma_{23}^*-\dfrac{\eta}{2}
\dfrac{\partial \sigma_{13}^*}{\partial \xi_1}
\bigg|_{\zeta_p=0}e^{-\zeta_q}+\eta^2\left(  
\dfrac{\partial^2 \sigma_{23}^*}{\partial \zeta_p^2} + \dfrac{1}{2}
\dfrac{\partial^2 \sigma_{13}^*}{\partial \xi_1 \partial\zeta_p }
\bigg|_{\zeta_p=0}e^{-\zeta_q}\right). 
\end{equation}
Re-writing the last formula in the original variables and substituting it into condition \eqref{eq_addcond} for $j=2$ we obtain at $x_2=0$
\begin{equation}\label{eq_bc_add1}
\sigma_{23}-a\left(\dfrac{\partial \sigma_{23}}{\partial x_2}  + \dfrac{\partial \sigma_{13}}{\partial x_1}\right)+a^2\left(   
\dfrac{\partial^2 \sigma_{23}}{\partial x_2^2}
+\dfrac{\partial^2 \sigma_{13}}{\partial x_1 \partial x_2}\right)=0.
\end{equation}
At the same time, using \eqref{eq_s32diff}, the original boundary condition \eqref{eq_bc} can be written as ($x_2=0$)
\begin{equation}\label{eq_bc_add2}
\sigma_{23}-\dfrac{a}{2} \dfrac{\partial \sigma_{13}}{\partial x_1}+a^2\left(   
\dfrac{\partial^2 \sigma_{23}}{\partial x_2^2}
+\dfrac{1}{2} \dfrac{\partial^2 \sigma_{13}}{\partial x_1 \partial x_2}\right)=f,
\end{equation}
which is in contradiction to the extra condition \eqref{eq_bc_add1}, following from \eqref{eq_addcond} at $j=2$. Thus, as might be expected, the asymptotic solutions, derived in Section 5, do not satisfy the integral formulation. Moreover, the violation of the condition \eqref{eq_addcond} at $j=2$ for the nonlocal stress $s_{23}$ means that the studied antiplane problem is ill-posed within the framework of the integral formulation. In fact, the reduction to a differential setup may be treated as a method for solving nonlocal equations initially presented in the integral form. Then, the conditions \eqref{eq_addcond} are related to the solvability of integral equations with the exponential  kernel $e^{-|x|/a}$, see e.g. \cite{polyanin2008handbook} for greater detail. The ill-poseddness of the static problem of a cantilever nonlocal elastic beam has been earlier revealed in \cite{romano2017constitutivebc}.

Nevertheless, the differential formulation analysed in this paper still has a visible potential. In particular, the spotted inconsistency will not be a feature of materials with only one of the stresses demonstrating nonlocal behaviour, analogously to the consideration in \cite{lurie2021dilatation}, dealing with simplified dilatational gradient elasticity.
In the latter case only one of the conditions \eqref{eq_addcond} has to be satisfied. Obviously, these observations hold true for more elaborated problems.

\section{Nonlocal shear surface wave}
Rewrite now formulae \eqref{eq_motfin} and \eqref{eq_bcfin} through displacement, neglecting the terms of $O(\eta^2)$. Then, in the homogeneous case $(f=0)$ 
\begin{equation}
\dfrac{\partial^2 u}{\partial x_1^2} + \dfrac{\partial^2 u}{\partial x_2^2}=\dfrac{1}{c_2^2}\dfrac{\partial^2 u}{\partial t^2},
\end{equation}
with
\begin{equation}
\left( \dfrac{\partial u}{\partial x_2} -\dfrac{a}{2} \dfrac{\partial^2 u}{\partial x_1^2} \right)\bigg|_{x_2=0}=0.
\end{equation}
Taking the displacement $u$ as
\begin{equation}
u = U e^{ik(x_1-ct)-k\gamma x_2},
\end{equation}
we obtain
\begin{equation}
\gamma = \sqrt{1-\dfrac{c^2}{c_2^2}},
\end{equation}
with the associated dispersion relation
\begin{equation}\label{eq_ROur}
C = \sqrt{1-\dfrac{1}{4}K^2},
\end{equation}
where, as above, 
\begin{equation}
C = \dfrac{c}{c_2} \quad \mbox{and} \quad K=ak.
\end{equation}
Over the long-wave region of interest $k \sim 1/l$.

We can also adapt for the considered case the refined boundary conditions derived in \cite{chebakov2016nonlocalhalfspace} within 3D integral nonlocal framework for a Gaussian kernel, giving
\begin{equation}
\left(\dfrac{\partial u}{\partial x_2}-\dfrac{a}{2\sqrt{\pi}} \dfrac{\partial^2 u}{\partial x_1^2} \right) \bigg|_{x_2=0}=0.
\end{equation}
This results in a similar dispersion relation
\begin{equation}\label{eq_RRoman}
C = \sqrt{1-\dfrac{1}{4\pi}K^2},
\end{equation}
which also supports a non-local shear surface wave.

\begin{figure}[H]
\centering
\begin{overpic}[scale=0.4,tics=5]{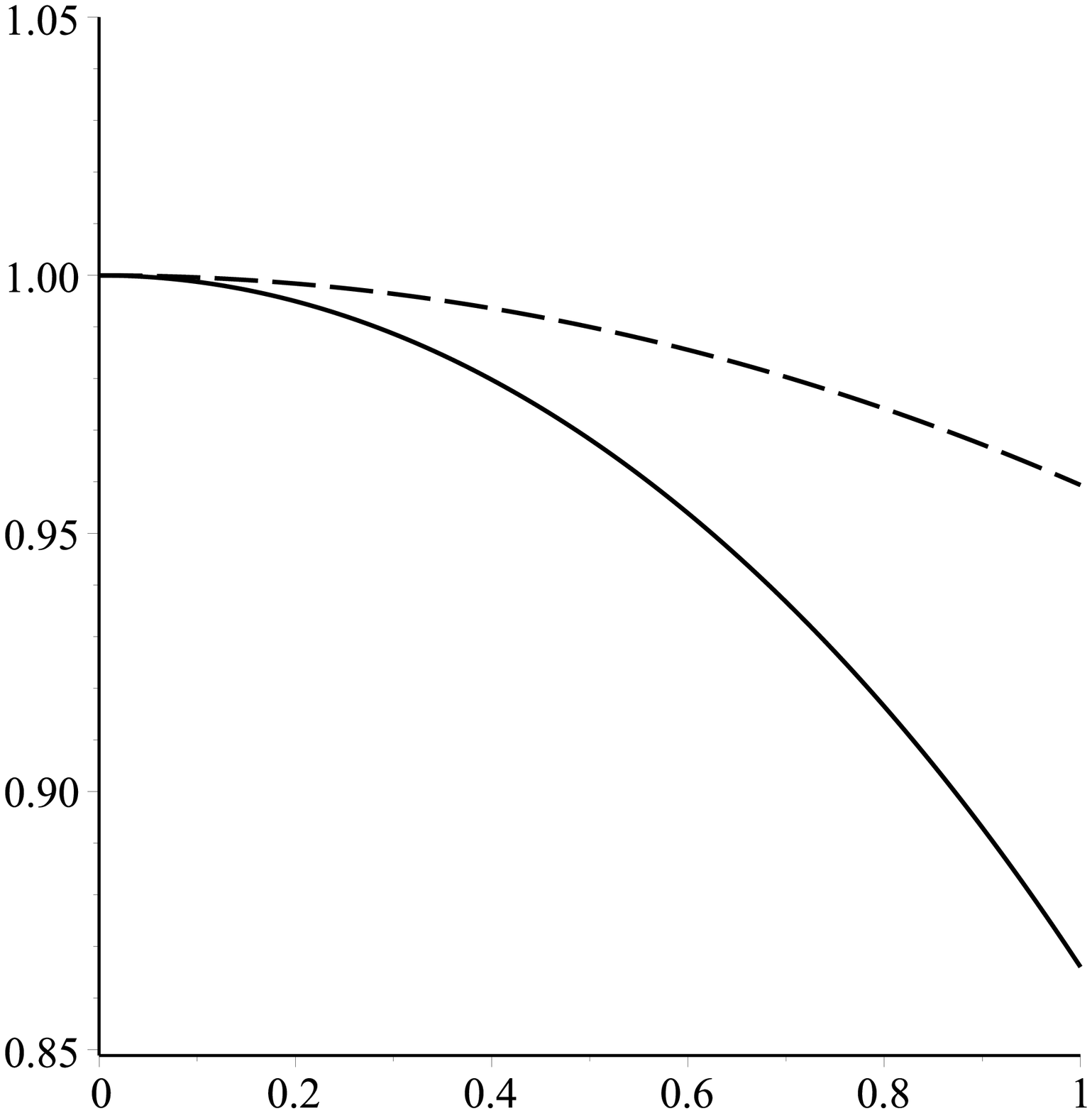}
  \put(-5,50){$C$}  
  \put(50,-5){$K$}  
\end{overpic}
  \medskip
\caption{Dispersion curves for the relations \eqref{eq_ROur} (solid line) and \eqref{eq_RRoman} (dashed line). }
\label{fig_speeds}
\end{figure}

The dispersion curves for relations \eqref{eq_ROur} and \eqref{eq_RRoman} are plotted in Figure \ref{fig_speeds}.
As expected, the dispersive effect due to nonlocality is visible at non-zero wave numbers.

\section{Conclusions}\label{conclusions}

The counter-example presented in Section 2 shows that the integral and differential models in nonlocal elasticity are not equivalent. This example is inspired by a related plane-strain time-harmonic problem, considered in \cite{eringen1983differential}. Undoubtedly, the solution in \cite{eringen1983differential} also does not  satisfy the equations of motion in nonlocal stresses. 

The full differential nonlocal formulation studied in Section 3 supports an antiplane surface wave, in contrast to the conventional differential setup, neglecting the effect of boundary layers associated with nonlocal stresses. The boundary layers in Section 3 are similar to those, determined from the nonlocal integral formulation with the Gaussian kernel \cite{chebakov2016nonlocalhalfspace, chebakov2017nonlocalplate}. 

Although boundary layers in nonlocal and gradient elasticity demonstrate a formal similarity, their physical nature appears to be quite different. In nonlocal elasticity the boundary layers are due to a rapid change in the integration domain near the surface, characteristic of localised kernels, involving a small internal size. On the contrary, singular perturbation in the differential equations of gradient elasticity usually come as higher order terms at homogenisation of periodic structures, resulting in the so-called spurious boundary layers, e.g. see \cite{pichugin2008asymptotic, pichugin2008gradient,
kaplunov2009rational}. The latter are the short-wave side products of long-wave asymptotic procedures. Spurious solutions also arise in refined theories for thin elastic plates and shells, e.g. see \cite{goldenveizer1993timoshenko}. The analogy between asymptotic schemes for thin and periodic structures is addressed in \cite{craster2014long}.

Another delicate issue is that the basic integral relations for nonlocal stresses generally need to be revisited when tackling boundary layers. Indeed, sticking with the assumption that the continuous integral relations originate from the homogenisation of the associated discrete chain or lattice, e.g. see \cite{eringen1977relation}, we may expect that the asymptotic procedure fails at the scale of an internal size, characteristic for nonlocal boundary layers. 
This stimulates the development of more elaborated approaches, combining discrete boundary layers and continuous outer solutions.
However, it might be expected that the effect will be only on coefficients at higher-order terms in effective boundary conditions.


%
\bibliographystyle{abbrv}

\bibliography{main}

\end{document}